\documentclass[a4paper,10pt]{article}
\mathcode`\<="4268     
\mathcode`\>="5269     
\mathcode`\:="603A     
\mathchardef\gt="313E  
\mathchardef\lt="313C  
\usepackage{mathbbol}
\usepackage{amsthm}
\usepackage{bussproofs}
\usepackage{tikz}
\theoremstyle{definition}

\def\textcdot{\raisebox{.5ex}{\texttt{.}}}
\def\spica{\mathrel{\texttt{))}}}
\def\spice{\mathrel{\texttt{)\textcdot(}}}
\def\spici{\mathrel{\texttt{()}}}
\def\spico{\mathrel{\texttt{(\textcdot(}}}

\def\aa#1#2{\ensuremath{\catebf{A}_{#1#2}}}
\def\ee#1#2{\ensuremath{\catebf{E}_{#1#2}}}
\def\ii#1#2{\ensuremath{\catebf{I}_{#1#2}}}
\def\oo#1#2{\ensuremath{\catebf{O}_{#1#2}}}

\def\AA#1#2{\ensuremath{\xymatrix{#1\ar[r]&#2}}}
\def\EE#1#2{\ensuremath{\xymatrix{#1\ar[r]&\b&#2\ar[l]}}}
\def\II#1#2{\ensuremath{\xymatrix{#1&\b\ar[l]\ar[r]&#2}}}
\def\OO#1#2{\ensuremath{\xymatrix{#1&\b\ar[l]\ar[r]&\b&#2\ar[l]}}}
\def\Ao#1#2{\ensuremath{\xymatrix{#2&#1\ar[l]}}}
\def\Eo#1#2{\ensuremath{\xymatrix{#2\ar[r]&\b&#1\ar[l]}}}
\def\Io#1#2{\ensuremath{\xymatrix{#2&\b\ar[l]\ar[r]&#1}}}
\def\Oo#1#2{\ensuremath{\xymatrix{#2\ar[r]&\b&\b\ar[l]\ar[r]&#1}}}
\newtheorem{theorem}{Theorem}[section]
\newtheorem{lemma}[theorem]{Lemma}
\newtheorem{proposition}[theorem]{Proposition}

\newtheorem{definition}[theorem]{Definition}
\newtheorem{remark}[theorem]{Remark}
\newtheorem{notation}[theorem]{Notation}
\newtheorem{example}[theorem]{Example}
\newtheorem{examples}[theorem]{Examples}
\newtheorem{terminology}[theorem]{Terminology}
\usepackage[english]{babel}
\usepackage[all]{xy}
\usepackage{txfonts}
\usepackage{dsfont}
\usepackage{bbm}
\usepackage{latexsym}
\title{A diagrammatic calculus of syllogisms\footnote{The final publication is available at springerlink.com http://link.springer.com/article/10.1007\%2Fs10849-011-9156-7\#page-1}}
\author{Ruggero Pagnan\\
DISI, University of Genova, Italy\\
\texttt{ruggero.pagnan@disi.unige.it}}
\date{}
\begin{document}
\maketitle
\begin{abstract}
A diagrammatic logical calculus for the syllogistic reasoning 
is introduced and discussed. We prove that a syllogism is valid if 
and only if it is provable in the calculus.
\end{abstract}
\def\a{\ensuremath{(a)}}
\def\adj#1#2{\ensuremath{\xymatrix@C.2cm{#1\ar@{-|}[r]&#2}}}
%
\def\adjpair#1#2#3#4#5{\ensuremath{\xymatrix{#1\ar@/^/[r]^{#3}_{\hole}="1"&#2\ar@/^/[l]^{#4}_{\hole}="2" \ar@{} "2";"1"|(.3){#5}}}}
%
\def\and{\ensuremath{\wedge}}
%
\def\arr{\ensuremath{\rightarrow}}
%
\def\Arr{\ensuremath{\Rightarrow}}
%
%
%
\def\arrow#1{\ensuremath{\cateb{#1}^{\arr}}}
%
\def\as{\ensuremath{\ast}}
%
%
\def\b{\ensuremath{\bullet}}
%
\def\BC{\ensuremath{(BC)}}
%
\def\BCd{\ensuremath{(BC)_d}}
\def\beps{\ensuremath{\backepsilon}}
%
\def\bicat#1{\ensuremath{\mathcal{#1}}}
%
\def\bf#1{\ensuremath{\mathbf{#1}}}
%
\def\bpararrow#1#2#3#4{\ensuremath{\xymatrix{#1\ar@<.9ex>[rr]^{#3}\ar@<-.9ex>[rr]_{#4}&&#2}}}
%
\def\bsquare#1#2#3#4#5#6#7#8{\xymatrix{#1\ar[dd]_{#8}\ar[rr]^#5&&#2\ar[dd]^{#6}\\
\\
#3\ar[rr]_{#7}&&#4}}
\def\bu{\ensuremath{\bullet}}
\def\btoind#1{#1\index{#1@\textbf{#1}}}
%
\def\card#1{\ensuremath{|#1|}}
\def\cate#1{\ensuremath{\mathcal{#1}}}
\def\cateb#1{\ensuremath{\mathbbm{#1}}}
%
\def\catebf#1{\ensuremath{\mathbf{#1}}}
\def\catepz#1{\ensuremath{\mathpzc{#1}}}
%
\def\cateop#1{\ensuremath{\cate#1^{op}}}
\def\comparrow#1#2#3#4#5#6{\ensuremath{\xymatrix{#6#1\ar[r]#4&#2\ar[r]#5&#3}}}
\def\cons#1{\ensuremath{\newdir{|>}{%
!/4.5pt/@{|}}\xymatrix@1@C=.3cm{\ar@{|=}[r]_{#1}&}}}
\def\D{\ensuremath{(D)}}
\def\enunciato#1#2#3#4#5#6#7{\newtheorem{#4}#5{#1}#6                                              
\begin{#4}#2\label{#7}
#3
\end{#4}}
%
\def\dfn#1#2#3#4#5#6#7{\newtheorem{#4}#5{#1}#6                                              
\begin{#4}#2\label{#7}
\emph{#3}
\end{#4}}
%
\def\epiarrow#1#2#3#4{\ensuremath{\newdir{|>}{%
!/4.5pt/@{|}*:(1,-.2)@^{>}*:(1,+.2)@_{>}}\xymatrix{#3#1\ar@{-|>}[r]#4&#2}}}
\def\epitip{\ensuremath{\newdir{|>}{%
!/4.5pt/@{|}*:(1,-.2)@^{>}*:(1,+.2)@_{>}}}}
\def\eps#1{\ensuremath{\epsilon#1}}
\def\esempio#1#2#3#4#5#6#7{\newtheorem{#4}#5{#1}#6                                              
\begin{#4}#2\label{#7}
\textup{#3}
\begin{flushright}
$\blacksquare$
\end{flushright}
\end{#4}}
%
\def\et#1{\ensuremath{\eta#1}}
\def\Ex#1{\mbox{\boldmath\ensuremath{\exists}}{#1}}
\def\ex#1{\ensuremath{\exists#1}}
\def\F{\ensuremath{(F)}}
\def\farrow#1#2#3{\ensuremath{\xymatrix{#3:#1\ar[r]&#2}}}
\def\Fi{\ensuremath{\Phi}}
%
\def\fib#1#2#3{\ensuremath{
\catebf{#3}:\cateb{#1}\arr\cateb{#2}}}
%
%
\def\fibs#1#2#3#4{\ensuremath{\begin{array}{ll}
\slice{\catebf{#1}}{#2}\\
\,\downarrow^{\catebf{#4}}\\
\catebf{#3}
\end{array}}}
\def\fibss#1#2#3#4#5{\ensuremath{\begin{array}{ll}
\slice{\catebf{#1}}{#2}\\
\,\downarrow^{#5}\\
\slice{\catebf{#3}}{#4}
\end{array}}}
\def\forev#1{\ensuremath{\forall~#1,~\forall~}}
\def\fract#1{\ensuremath{\cateb{#1}[\Sigma^{-1}]}}
%
\def\G#1{\ensuremath{\int{#1}}}
%
%
\def\harrow{\ensuremath{\rightharpoonup}}
%
\def\hook{\ensuremath{\hookrightarrow}}
\def\hset#1#2#3{\ensuremath{\cate{#1}(#2,#3)}} 
\def\Im#1{\ensuremath{\mathit{Im}#1}}
\def\implies{\ensuremath{\Rightarrow}}
\def\indprod#1#2{\ensuremath{\prod#1#2}}
%
\def\infer#1#2{\ensuremath{\left.\begin{array}{cc}
#1\\
\xymatrix{\ar@{-}@<-.3ex>[rrrr]\ar@{-}@<.3ex>[rrrr]&&&&}\\
#2\\
\end{array}\right.}}
%
%
\def\it#1{\ensuremath{\mathit{#1}}}
%
\def\l#1{\ensuremath{#1_{\bot}}}
\def\larrow{\ensuremath{\leftarrow}}
\def\lcomparrow#1#2#3#4#5{\ensuremath{\xymatrix{#1\ar[rr]#4&&#2\ar[rr]#5&&#3}}}
\def\lfract#1{\ensuremath{[\Sigma^{-1}]\cateb{#1}}}
%
\def\lhookarrow#1#2#3{\ensuremath{\xymatrix{#1\ar@{^{(}->}[rr]^{#3}&&#2}}}
\def\mod#1#2{\ensuremath{\cateb{Mod}(#1,#2)}}
\def\Mono#1{\ensuremath{\mathit{Mono}#1}}
\def\monoarrow#1#2#3#4{\ensuremath{\newdir{
>}{{}*!/-5pt/@{>}}\xymatrix{#3#1\ar@{ >->}[r]#4&#2}}}
\def\monotail{\ensuremath{\newdir{
>}{{}*!/-5pt/@{>}}}}
\def\N{\ensuremath{(N)}}
\def\natarrow#1#2#3{\ensuremath{\newdir{|>}{%
!/4.5pt/@{|}*:(1,-.2)@^{>}*:(1,+.2)@_{>}}\xymatrix{#3:#1\ar@{=|>}[r]&#2}}}
\def\natbody{\ensuremath{\newdir{|>}{%
!/4.5pt/@{|}*:(1,-.2)@^{>}*:(1,+.2)@_{>}}}}
\makeatletter
\def\natpararrow#1#2#3#4#5{\@ifnextchar(
 {\Natpararrow{#1}{#2}{#3}{#4}{#5}}
 {\Natpararrow{#1}{#2}{#3}{#4}{#5}(.5)}}
\makeatother
\def\Natpararrow#1#2#3#4#5(#6){\ensuremath{\newdir{|>}{%
!/4.5pt/@{|}*:(1,-.2)@^{>}*:(1,+.2)@_{>}}
\xymatrix{#1\ar@<1.5ex>[rr]^(#6){#3}|(#6){\vrule height0pt depth.7ex
width0pt}="a"\ar@<-1.5ex>[rr]_(#6){#4}|(#6){\vrule height1.5ex width0pt}="b"&&#2
\ar@{=|>} "a";"b"#5}}}
\def\om{\ensuremath{\omega}}
%
\def\ov#1{\ensuremath{\overline{#1}}}
%
\def\ovarr#1{\ensuremath{\overrightarrow{#1}}}
%
\def\pararrow#1#2#3#4{\ensuremath{\xymatrix{#1\ar@<.8ex>[r]^{#3}\ar@<-.8ex>[r]_{#4}&#2}}}
\def\p#1{\ensuremath{\mathds{P}_{\catepz{S}}}(#1)}
\def\P#1{\ensuremath{\mathnormal{P}(#1)}}
\def\pr#1{\textbf{Proof:~}#1 
\begin{flushright}
$\Box$
\end{flushright}
}
%
%
%
%
\def\rel#1#2#3{\ensuremath{\xymatrix{#1:#2\ar[r]
|-{\SelectTips{cm}{}\object@{|}} &#3}}}
%
%
\def\res#1#2{\ensuremath{#1_{\rbag #2}}}
\def\S#1{\ensuremath{\catepz{S}_{\cateb{#1}}}}
\def\sectoind#1#2{#1\index{#2!#1}}
%
\def\sinfer#1#2{\ensuremath{\left.\begin{array}{cc}
#1\\
\xymatrix{\ar@{-}[rrrr]&&&&}\\
#2\\
\end{array}\right.}}
%
\def\slice#1#2{\ensuremath{#1/{\textstyle#2}}}
\def\sslice#1#2{\ensuremath{#1//{\textstyle#2}}}
\def\square#1#2#3#4#5#6#7#8{$$\xymatrix{#1\ar[d]_{#8}\ar[r]^{#5}&#2\ar[d]^{#6}\\
#3\ar[r]_{#7}&#4}$$}
\def\sub#1#2#3{\ensuremath{#1(#2#3)~}}
\def\summ#1{\ensuremath{\sum#1}}
%
\def\ti#1{\ensuremath{\tilde#1}}
%
\def\tens{\ensuremath{\otimes}}
%
\def\teo#1{\ensuremath{\tau#1}}
\def\ter#1{\ensuremath{\tau_{\cate{#1}}}}
\def\teta{\ensuremath{\theta}}
%
%
\def\toind#1{#1\index{#1}}
%
\def\un#1{\underline{#1}}
\def\vcell{\ensuremath{\ar@<-1ex>@{}[r]_{\hole}="a"\ar@<+1ex>@{}[r]^{\hole}="b"}}
%
\def\veps{\ensuremath{\varepsilon}}
%
\def\vfi{\ensuremath{\varphi}}
\def\w{\ensuremath{\wedge}}
\def\wti#1{\ensuremath{\widetilde#1}}
\def\y#1{\ensuremath{y#1}}

\section{Introduction}
The aim of the present paper is that of introducing and discussing a
diagrammatic logical calculus for syllogistic reasoning. We present
suitable linear diagrammatic representations of the fundamental
Aristotelian categorical propositions and show that 
they are closed under the syllogistic canon of
inference which is the deletion of the middle term, so
implemented to let the formalism incorporate simultaneously a graphical 
appearance and a naive algorithmic nature, namely that no
specific knowledge or particular ability is needed in order to
understand it and use it.\\ 
Since our investigation is directed toward 
a formal approach to logical reasoning with diagrams we introduce a formal system SYLL for such a calculus. We prove that a syllogism is valid if and only if
it is provable in SYLL, so that in this sense the calculus is sound and complete. A similar result holds also for the syllogisms that are valid under existential import.
Because of the peculiar form of its diagrammatic syntax, the calculus 
supports a criterion for the rejection of the invalid syllogistic arguments on the base of which the easy retrieving of the traditional rules of the syllogism is possible. Moreover, we show that
the laws of the square of opposition are provable in SYLL.\\ 
In section~\ref{sec2} we introduce the basics of the syllogistics and describe the diagrammatic logical system VENN based on the Venn-Peirce diagrams, from~\cite{MR1271699} and~\cite{MR1408440}.\\ 
In section~\ref{sec4} we prove the  previous claims about soundness and completeness, while
the laws of the square of opposition are discussed in section~\ref{sec5}.\\ 
The possibility of extending the calculus to $n$-term syllogistic inferences is briefly discussed in section~\ref{disc}. We point out that other linear diagrammatic formalisms for the syllogistic reasoning exist, notably~\cite{Smyth},~\cite{MR1149957},~\cite{MR1776228} and 
that a category-theoretic point of view is pursued in~\cite{MR1271697}. 
We are aware of possible, interesting directions of investigation in connection with Peirce's Existential Graphs, see~\cite{MR0497852}, while further directions of investigation could be pursued in connection with computer science, see~\cite{Rayside:2001:SSO:381473.381485}.\\
\linebreak
I acknowledge Pino Rosolini for the many useful conversations and the anonymous referees for their many
valuable comments.

\section{Syllogistic reasoning with diagrams}\label{sec2}
Syllogistics in its original form dates back to Aristotle, who formalized it as a logical system in the Prior Analytics. A fundamental reference on the subject is~\cite{Lukasiewicz}, for example. We here recall some of the basics of syllogistics in its traditional, medieval systematization. In doing this we introduce some of the terminology and notations that will be useful in the sequel of the paper. Moreover, we recall the fundamentals of a formal diagrammatic-theoretic approach to syllogistics by reviewing briefly
the formal system VENN, see~\cite{MR1271699},~\cite{MR1408440} and~\cite{MR1312613}.\\

We refer to nouns, adjectives or more
complicated meaningful expressions of the natural language as \emph{terms},
and denote them by using upper case latin letters which we refer to as
\emph{term-variables}. Since Aristotle, the following four schemes of propositions were
recognized as fundamental throughout the research in logic:
\[\begin{tabular}{ll}
\aa{A}{B}: Each $A$ is $B$ & universal affirmative proposition\\
\\
\ee{A}{B}: No $A$ is $B$ & universal negative proposition\\
\\
\ii{A}{B}: Some $A$ is $B$ & particular affirmative proposition\\
\\
\oo{A}{B}: Some $A$ is not $B$ & particular negative proposition
\end{tabular}\]
Following the tradition, henceforth we refer to them as
\emph{categorical propositions}. In each of them, the term-variable $A$ is
the \emph{subject} whereas the term-variable $B$ is the \emph{predicate} of
the proposition. Thus ``Each dog is black'' , ``No cat is
white'' , ``Each baby that cry is polite'' are examples of
categorical propositions.\\ 
A \emph{syllogism} is a logical consequence
that involves three categorical propositions that are distinguished in
\emph{first premise}, \emph{second premise} and \emph{conclusion}.
Moreover, a syllogism involves three term-variables $S$, $P$ and
$M$ in the following precise way: $M$ does not
occur in the conclusion whereas, according to the traditional way of
writing syllogisms, $P$ occurs in the first premise and 
$S$ occurs in the second premise. The term-variables $S$ and $P$ 
occur as the subject and predicate of the conclusion, respectively,
and are also referred to as \emph{minor term} and \emph{major term} of
the syllogism, whereas $M$ is also referred to as \emph{middle
term}.
\begin{remark}
What we are simply referring to as syllogisms are
\emph{traditional syllogisms} in the terminology
of~\cite{Lukasiewicz}, where a detailed discussion of the
difference between this notion and that of \emph{Aristotelian syllogism}
can be found. Such a difference will not affect the present
treatment. We here only mention
that in strict terms an Aristotelian syllogism is a proposition
of the type ``If A and B, then C'', whereas a traditional syllogism
is a logical consequence with two premises and one conclusion like
``A, B therefore C'', which in its entirety does not form a
compound proposition. An Aristotelian syllogism can
either be true or false whereas a traditional syllogism can either be
valid or not, in the sense of Tarski, see~\cite{mendelson.iml}.
\end{remark}
The \emph{mood} of a syllogism is the sequence of the kinds of the
categorical propositions by which it is formed. The \emph{figure} of a
syllogism is the position of the term-variables $S$, $P$ and $M$ in it. There
are four possible figures as shown in the following table:
\begin{eqnarray}\label{figs}
\begin{tabular}{|l|c|c|c|c|}
\hline
& fig. 1 & fig. 2 & fig. 3 & fig. 4\\
\hline
first premise & MP & PM & MP & PM\\
\hline
second premise & SM & SM & MS & MS\\
\hline
conclusion & SP & SP & SP & SP\\
\hline
\end{tabular}
\end{eqnarray}
A syllogism is completely determined by its mood and by its figure
together. We write syllogisms so that their mood and figure can be
promptly retrieved, by also letting the symbol $\vdash$ separate the
premises from the conclusion.
For example, in the syllogism
\[\catebf{A}_{MP},\catebf{A}_{SM}\vdash\catebf{A}_{S P}\]
one can recognize from left to right the first premise, the
second premise and the conclusion, its mood, which is
\catebf{AAA}, and its figure, which is the first. 
The following tables list the syllogisms that are known to be valid since Aristotle.
\begin{eqnarray}\label{questa}{\textrm{\footnotesize
\begin{tabular}{|l|l|l|l|l|}
\hline
Fig. 1 & Fig. 2 & Fig. 3 & Fig. 4 \\
\hline
$\catebf{A}_{MP},\catebf{A}_{SM}\vdash\catebf{A}_{SP}$&$\catebf{E}_{PM},\catebf{A}_{SM}\vdash
\catebf{E}_{SP}$
&$\catebf{I}_{MP},\catebf{A}_{MS}\vdash\catebf{I}_{SP}$&$\catebf{A}_{PM},\catebf{E}_{MS}\vdash
\catebf{E}_{SP}$\\
$\catebf{E}_{MP},\catebf{A}_{SM}\vdash\catebf{E}_{SP}$&$\catebf{A}_{PM},\catebf{E}_{SM}\vdash
\catebf{E}_{SP}$
&$\catebf{A}_{MP},\catebf{I}_{MS}\vdash\catebf{I}_{SP}$&
$\catebf{I}_{PM},\catebf{A}_{MS}\vdash\catebf{I}_{SP}$\\
$\catebf{A}_{MP},\catebf{I}_{SM}\vdash\catebf{I}_{SP}$&$\catebf{E}_{PM},\catebf{I}_{SM}\vdash
\catebf{O}_{SP}$
&$\catebf{O}_{MP},\catebf{A}_{MS}\vdash\catebf{O}_{SP}$&
$\catebf{E}_{PM},\catebf{I}_{MS}\vdash\catebf{O}_{SP}$\\
$\catebf{E}_{MP},\catebf{I}_{SM}\vdash\catebf{O}_{SP}$&$\catebf{A}_{PM},\catebf{O}_{SM}\vdash
\catebf{O}_{SP}$
&$\catebf{E}_{MP},\catebf{I}_{MS}\vdash\catebf{O}_{SP}$&\\
\hline
\end{tabular}}}
\end{eqnarray}
\begin{eqnarray}\label{equesta}{\textrm{\footnotesize
\begin{tabular}{|l|l|l|l|l|}
\hline
Fig. 1 & Fig. 2 & Fig. 3 & Fig. 4 & assumption\\
\hline
$\aa{M}{P},\aa{S}{M}\vdash\ii{S}{P}$&$\aa{P}{M},\ee{S}{M}\vdash\oo{S}{P}$&
&$\aa{P}{M},\ee{M}{S}\vdash\oo{S}{P}$& some S exists\\
$\ee{M}{P},\aa{S}{M}\vdash\oo{S}{P}$&$\ee{P}{M},\aa{S}{M}\vdash\oo{S}{P}$&&&
some S exists\\
&&$\aa{M}{P},\aa{M}{S}\vdash\ii{S}{P}$&$\ee{P}{M},\aa{M}{S}\vdash\oo{S}{P}$&
some M exists\\
&&$\ee{M}{P},\aa{M}{S}\vdash\oo{S}{P}$&& some M exists\\
&&&$\aa{P}{M},\aa{M}{S}\vdash\ii{S}{P}$& some P exists\\
\hline
\end{tabular}}}
\end{eqnarray}
They are 24 in total, divided into two groups of 15 and of 9. 
Those syllogisms in the second group are also said to be \emph{strengthened}, 
or valid under \emph{existential import}, 
which is an explicit assumption of existence of some $S$,
$M$ or $P$, as indicated.
\begin{terminology}
The syllogisms in table~(\ref{questa}) will be henceforth referred to simply as syllogisms, those in table~(\ref{equesta}) will be referred to as strengthened syllogisms.
\end{terminology}

The study of logical reasoning requires to understand that the valid reasoning consists in the correct manipulation of the information no matter the nature of the symbolic medium, being it diagrammatic, linguistic or both, that is heterogeneous. The correct manipulation of the information is supported by the employment of sound rules of inference within suitable formal logical systems. In pursuing a diagrammatic approach to logical reasoning one point is that of making formal what is usually considered as heuristic.\\ 
A paradigmatic example of a formal diagrammatic logical system supporting syllogistic reasoning is the system VENN 
in~\cite{MR1408440} and~\cite{MR1271699}, which is based on the Venn-Peirce diagrams, see also~\cite{MR1408441}. Here we briefly describe VENN by recalling its diagrammatic syntactic primitives and its rules of inference. For further details we refer the reader to loc. cit.

\begin{definition}
The \emph{diagrammatic primitives} of VENN are the following distinct syntactic objects:
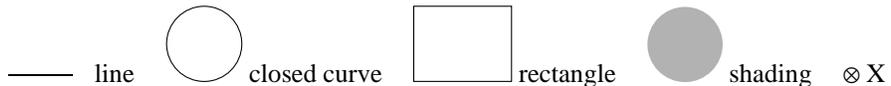
\[\begin{tabular}{ccccc}
$\xymatrix{\ar@{-}[r]&}$ line &
$\begin{tikzpicture}
\draw (0,0) circle (.5);\end{tikzpicture}$ closed curve &
$\begin{tikzpicture}
\draw (0,0) rectangle (1.3,1);
\end{tikzpicture}$ rectangle &
$\begin{tikzpicture}[fill=black!30]
\fill (0,0) circle (.5); 
\end{tikzpicture}$ shading &
$\otimes$ X
\end{tabular}\]
A \emph{diagram} of VENN is any finite combination of diagrammatic primitives. In particular, an X-\emph{sequence} is a diagram of alternating X's and lines with an X in each extremal position, e.g. $\xymatrix{\otimes\ar@{-}[r]&\otimes\ar@{-}[r]&\otimes}$.
Diagrams of VENN will be denoted by calligraphic upper case letters. For convenience, diagrams which are closed curves or rectangles will be denoted by upper case latin letters.
A \emph{region} is any enclosed area in a diagram. A \emph{basic region} is a region enclosed by a rectangle or a closed curve. A \emph{minimal region} is a region within which no other region is enclosed. 
\end{definition}
Regions of diagrams are meant to represent sets of individuals. In particular, a background rectangle is meant to represent a suitable universe of discourse.
A shaded region is understood as empty, whereas a region that contains an X-sequence is understood as non-empty. An X-sequence represents disjunctive existential statements. \\

The informal procedure for the verification of the validity of a syllogism through the employment of Venn-Peirce diagrams consists in drawing the diagrams for the premises as well as for the conclusion of the syllogism and see if it is possible to read off the latter from the former. If so, then the syllogism is valid, otherwise it is not. The crucial step consists in understanding if the diagram for the conclusion is ``contained'' in the diagrams for the premises or not. Making formal such a procedure requires to take into account how each diagram has been constructed, to understand for which reasons the construction of those diagrams is permitted and in particular to understand the reasons why the derivation of the diagram for the conclusion is permitted. We skip the formal presentation of how to construct a well-formed diagram in VENN and we refer the reader to~\cite{MR1408440}.

\begin{definition}
The \emph{rules of inference} of VENN are the following: 
\begin{description}
\item[setup:] a well-formed diagram with no shadings or X-sequences can be asserted at any step of a proof.
\item[erasure:] a well-formed diagram \cate{E} is obtained from a well-formed diagram \cate{D} by this rule if \cate{E} results from either erasing a closed curve of \cate{D}, or a shading of some region of \cate{D}, or an entire X-sequence of \cate{D}. In the first case, a shading filling part of a minimal region must also be erased.
\item[extension of a sequence:] \cate{E} is obtained from \cate{D} by this rule if extra links have been added to some X-sequence of \cate{D}.
\item[erasure of links:] \cate{E} is obtained from \cate{D} by this rule if it results from \cate{D} by the erasure of links in some X-sequence that fall in shaded regions, provided that the remaining X's are reconnected.
\item[unification:] \cate{D} is obtained from $\cate{D}_1$ and $\cate{D}_2$ by this rule if 
every region of \cate{D} is in counterpart relation with a region of either $\cate{D}_1$ or $\cate{D}_2$ and conversely. If any region of \cate{D} is shaded or has an X-sequence, then it has a counterpart in either $\cate{D}_1$ or $\cate{D}_2$ which is also shaded or has an X-sequence and conversely.
\item[non-emptyness] \cate{D} is obtained from \cate{E} by this rule if it has been obtained by the addition of an X-sequence some link of which falls into every minimal region of \cate{E}.
\end{description}
\end{definition}

The following are the well-formed diagrams of VENN that correspond to the categorical proposition:

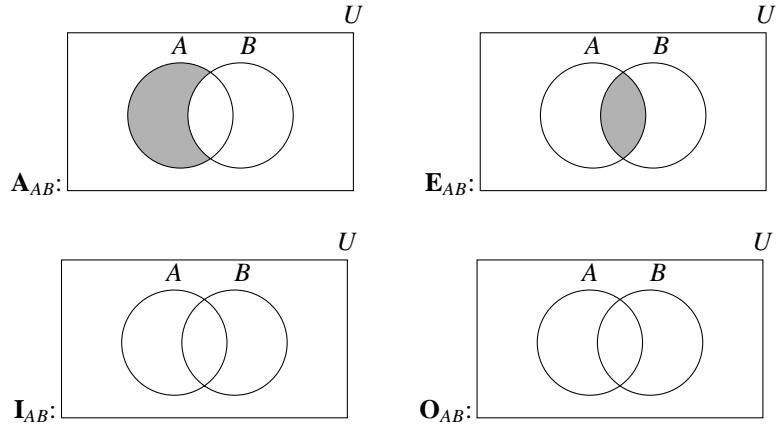
\[\catebf{A}_{AB}:
\begin{tikzpicture}[fill=black!30]
\scope
\clip (-1,-.5) rectangle (2.8,1.6)
      (1.3,.5) circle (.7);
\fill (.5,.5) circle (.7);
\endscope
\draw (.5,.5) circle (.7) (.5,1.2)  node [text=black,,above] {$A$}
      (1.3,.5) circle (.7) (1.4,1.2)  node [text=black,above] {$B$}
      (-1,-.5) rectangle (2.8,1.6) node [text=black,above] {$U$};
\end{tikzpicture}
\qquad
\catebf{E}_{AB}:
\begin{tikzpicture}[fill=black!30]
\scope
\clip 
      (.5,.5) circle (.7);
\fill (1.3,.5) circle (.7);
\endscope
\draw (.5,.5) circle (.7) (.5,1.2)  node [text=black,,above] {$A$}
      (1.3,.5) circle (.7) (1.4,1.2)  node [text=black,above] {$B$}
      (-1,-.5) rectangle (2.8,1.6) node [text=black,above] {$U$};
\end{tikzpicture}\]
\[\catebf{I}_{AB}:
\begin{tikzpicture}
\draw (.5,.5) circle (.7) (.5,1.2)  node [text=black,,above] {$A$}
      (1.3,.5) circle (.7) (1.4,1.2)  node [text=black,above] {$B$}
      (-1,-.5) rectangle (2.8,1.6) node [text=black,above] {$U$};
\put (23,12) {$\otimes$};
\end{tikzpicture}
\qquad
\catebf{O}_{AB}:
\begin{tikzpicture}
\draw (.5,.5) circle (.7) (.5,1.2)  node [text=black,,above] {$A$}
      (1.3,.5) circle (.7) (1.4,1.2)  node [text=black,above] {$B$}
      (-1,-.5) rectangle (2.8,1.6) node [text=black,above] {$U$};
\put (7,12) {$\otimes$};
\end{tikzpicture}\]

Here are two examples of formal proofs in VENN
in which rectangles have been omitted.

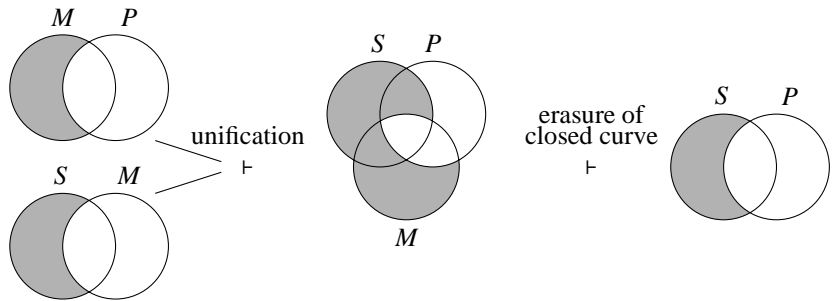
\def\tta[#1,#2]#3#4#5{\begin{scope}[xshift=#1,yshift=#2]
#5;
\draw (0,0) circle (20) (0,20) node [above] {$#3$}
(20,0) circle (20) (25,20) node [above] {$#4$};
\end{scope}}
\def\tA[#1,#2]#3#4;{\tta[#1,#2]{#3}{#4}{\begin{scope}
\clip (-30,-30) rectangle (50,45)
(20,0) circle (20);
\fill (0,0) circle (20);
\end{scope}}}
\def\tI[#1,#2]#3#4;{\tta[#1,#2]{#3}{#4}{\node at (10,0) {$\otimes$}}}

\def\tb(#1)#2#3{\begin{scope}[xshift=#1,yshift=0]
#2
\begin{scope}
\clip (-30,-30) rectangle (50,45)
(20,20) circle (20);
\fill (10,0) circle (20);
\end{scope}
\draw (0,20) circle (20) (0,40) node [above] {$S$}
 (20,20) circle (20) (20,40) node [above] {$P$}
 (10,0) circle (20) (10, -20) node [below] {$M$};
#3
\end{scope}}

\def\mt(#1)#2#3{\node at (#1,20) {#2};
\node at (#1,12) {#3};
\node at (#1,0) {$\vdash$}}
\def\mo(#1)#2{\node at (#1,12) {#2};
\node at (#1,0) {$\vdash$}}

\begin{itemize}
\item[-] $\aa{M}{P},\aa{S}{M}\vdash\aa{S}{P}$
\begin{center}
\begin{tikzpicture}[x=1pt,y=1pt,fill=black!30] 
\tA[0,30]MP;
\tA[0,-30]SM;
\draw (35,10) -- (60,2) ;
\draw (35,-10) -- (60,-2) ;
\mo(70){unification};
\tb(120){\begin{scope}
\clip (-30,-30) rectangle (50,45)
(10,0) circle (20);
\fill (0,20) circle (20);
\end{scope}}{}
\mt(200){erasure of}{closed curve};
\tA[250,0]SP;
\end{tikzpicture}
\end{center}

\item[-] $\aa{M}{P},\ii{S}{M}\vdash\ii{S}{P}$

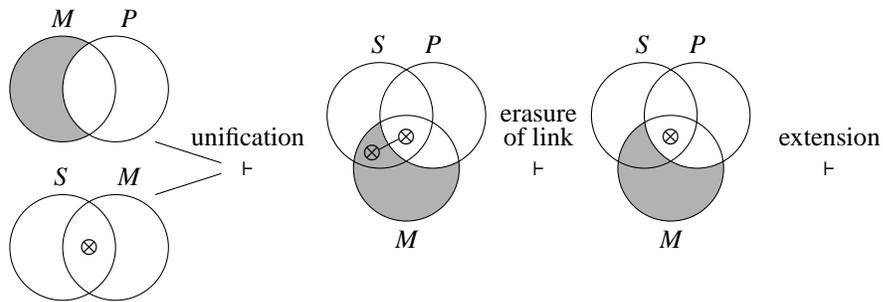
\begin{center}
\begin{tikzpicture}[x=1pt,y=1pt,fill=black!30]
\tA[0,30]MP;
\tI[0,-30]SM;
\draw (35,10) -- (60,2) ;
\draw (35,-10) -- (60,-2) ;
\mo(70){unification};
\tb(120){}{\node at (-3,6) {$\otimes$} ;
\node at (10,12) {$\otimes$} ;
\draw (-.5,7) -- (7,11) ;}
\mt(180){erasure}{of link};
\tb(220){}{\node at (10,12) {$\otimes$} ;}
\mo(290){extension};
\end{tikzpicture}

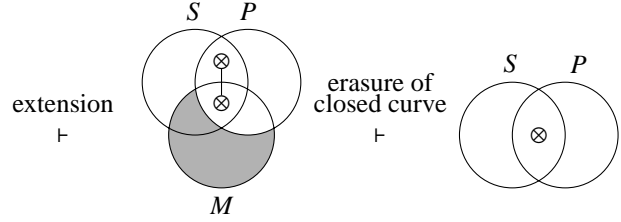
\begin{tikzpicture}[x=1pt,y=1pt,fill=black!30]
\mo(0){extension};
\tb(50){}{\node at (10,28) {$\otimes$} ;
\node at (10,12) {$\otimes$} ;
\draw (10,15) -- (10,25) ;}
\mt(120){erasure of}{closed curve};
\tI[170,0]SP;
\end{tikzpicture}
\end{center}

\end{itemize}

For the proof that VENN is sound and complete, we refer the reader to~\cite{MR1408440} and~\cite{MR1271699}.

\section{The calculus}\label{sec4}
In this section we introduce the formal system SYLL, 
supporting a diagrammatic logical calculus for the syllogistic reasoning. A feature of SYLL is that it is heterogeneous, in the sense that it consists of diagrammatic and linguistic syntactic objects together. We will prove that the calculus at issue is sound and complete, in the sense that a syllogism is valid if and only if it is provable in SYLL.

\begin{definition}\label{synt}
The \emph{diagrammatic primitives} of SYLL are the symbols $\arr$, $\leftarrow$, $\b$. The \emph{linguistic primitives} of SYLL consist of countably many term-variables $A,B,C,\ldots$. The \emph{syntactic primitives} of SYLL are the diagrammatic or linguistic primitives. To each scheme of categorical proposition we associate the following schemes of \emph{syllogistic diagrams}
\[\begin{tabular}{lllll}
\aa{A}{B}: & \AA{A}{B} &&
\ee{A}{B}: & \EE{A}{B}\\
\\
\ii{A}{B}: & \II{A}{B} &&
\oo{A}{B}: & \OO{A}{B}
\end{tabular}\]
to be read analogously. A \emph{diagram} of SYLL is a finite list of arrow symbols separated by a single bullet symbol or term-variable, beginning and ending at a term-variable. The \emph{reversal} of a given diagram is the diagram obtained by specular symmetry. A \emph{part} of a diagram 
is a finite list of consecutive components of a diagram.
\end{definition}

\begin{examples}
The lists $A$, $A\arr X$, $A\leftarrow A$, $A\arr\b\arr B$, $X\arr Y\arr\b\leftarrow X$ are examples of diagrams. Their reversals are $A$, $X\leftarrow A$, $A\arr A$, $B\leftarrow\b\leftarrow A$ and $X\arr\b\leftarrow Y\leftarrow X$, respectively.
The reversals of the syllogistic diagrams are the diagrams
\[\begin{tabular}{llll}
\Ao{A}{B} && \Eo{A}{B}\\
\\
\Io{A}{B} && \Oo{A}{B}
\end{tabular}\]
Every diagram is a part of itself. In general, a part of a diagram need not be a diagram, e.g. $A\arr$ is a part of the diagram $A\arr\b\leftarrow B$ and it is not a diagram, since it does not end at a term-variable.
\end{examples}

\begin{notation}
Parts of diagrams will be henceforth denoted by calligraphic upper case letters such as $\cate{D},\cate{E}$, etc. In order to distinguish explicitly a part with respect to a whole diagram, we adopt a heterogeneous notation mixing calligraphic upper case letters and syntactic primitives. For example, the writing $\cate{D}\arr A$ refers to a diagram in which the part $\arr A$ has been distinguished with respect to the remaining part \cate{D}. Thus, it may be the case that the whole diagram looks like $X\leftarrow\b\arr A$, $S\arr A$ or $B\arr \b\arr A$ for example, so that the part \cate{D} would be $X\leftarrow \b$, $S$, $B\arr \b$, respectively.
\end{notation}

\begin{definition}
A \emph{concatenable pair} of diagrams is a pair of diagrams $(\cate{D}A,A\cate{E})$ or $(A\cate{E},\cate{D}A)$ whose \emph{concatenation} is, in both cases, the diagram $\cate{D}A\cate{E}$ which is obtained by overlapping its components on the common extremal term-variable $A$. A \emph{composable pair} of diagrams is a concatenable pair $(\cate{D}\arr A,A\arr\cate{E})$ or $(A\arr\cate{E},\cate{D}\arr A)$, $(\cate{D}\leftarrow A,A\leftarrow\cate{E})$ or $(A\leftarrow\cate{E},\cate{D}\leftarrow A)$. In the first two cases, a composable pair gives rise to a \emph{composite}
$\cate{D}\arr\cate{E}$ obtained by substituting the part $\arr A\arr$ in the concatenation $\cate{D}\arr A\arr \cate{E}$ with the sole, accordingly oriented, arrow symbol $\arr$.
Analogously, in the second two cases, a composable pair gives rise to a composite diagram $\cate{D}\leftarrow\cate{E}$. For every natural number $n$, $n\geq 3$, a \emph{concatenable $n$-tuple} is an $n$-tuple of diagrams $(\cate{E}_1,\cate{E}_2,\ldots,\cate{E}_n)$ in which, for every $1\leq i\lt n$, the pairs $(\cate{E}_i,\cate{E}_{i+1})$ are concatenable pairs of the same form 
$(\cate{D}A,A\cate{E})$ or $(A\cate{E},\cate{D}A)$. A \emph{composable $n$-tuple} is a concatenable $n$-tuple of diagrams $(\cate{E}_1,\cate{E}_2,\ldots,\cate{E}_n)$ in which, for every $1\leq i\lt n$,  $(\cate{E}_i,\cate{E}_{i+1})$ is a composable pair. Composition of diagrams extends to composable $n$-tuples through the calculation of pairwise composites.
\end{definition}

\begin{examples}
For every term-variable $A$, $(A,A)$ is a concatenable pair whose concatenation is the diagram $A$. It is not a composable pair since no arrow symbols occur. The pair $(A\leftarrow B, X\arr B)$ is not concatenable, thus not composable, whereas the pair $(A\leftarrow B,B\leftarrow X)$ is concatenable and composable, with composite $A\leftarrow X$. The pair $(X\arr A, A\leftarrow B)$ is concatenable but not composable.
The pair $(X\leftarrow B,B\leftarrow X)$ is concatenable in two different ways by overlapping its components either on $B$ or on $X$. Also, it is composable in two different ways giving rise to either the composite $X\leftarrow X$ or the composite $B\leftarrow B$, respectively. The $3$-tuple $(X\leftarrow B, B\leftarrow X, A\arr B)$ is concatenable to $A\arr B\leftarrow X\leftarrow B$ but it is not composable since the pair $(B\leftarrow X,A\arr B)$ is not composable. The $3$-tuple $(X\leftarrow A, B\arr X, X\arr A)$ is not concatenable, since the pair $(X\leftarrow A, B\arr X)$ is concatenable by overlapping its components on $X$ in extremal ``external'' position, whereas the pair $(B\arr X, X\arr A)$ is concatenable by overlapping its components on $X$ in extremal ``internal'' position.
The $3$-tuple 
$(A\leftarrow \b\arr X,X\arr\b\leftarrow T,T\leftarrow\b\arr H)$
is concatenable to 
$A\leftarrow \b\arr X\arr\b\leftarrow T\leftarrow \b\arr H$
and composable to 
$A\leftarrow\b\arr X\arr\b\leftarrow\b\arr H$.
\end{examples}

\begin{definition}
A \emph{well-formed diagram} of SYLL is defined inductively as follows:
\begin{itemize}
\item[(i)] a syllogistic diagram is a well-formed diagram.
\item[(ii)] the reversal of a syllogistic diagram is a well-formed diagram.
\item[(iii)] a diagram which is the concatenation of a concatenable pair whose components are 
well-formed diagrams is a well-formed diagram.
\end{itemize}
\end{definition}

\begin{remark}
Well-formed diagrams are not closed under composition. Indeed, it suffices to consider the composable pair $(X\arr\b\leftarrow\b\arr A,A\arr\b\leftarrow\b\arr Y)$
for example, whose components are well-formed but give rise to the composite diagram
$X\arr\b\leftarrow\b\arr\b\leftarrow\b\arr Y$
which is not well-formed.
\end{remark}

The intuition about how to use the syllogistic diagrams and their reversals to verify the validity of syllogisms is that, given a syllogism, one considers the three syllogistic diagrams or reversals to represent the first premise, the second premise and the conclusion of the syllogism. These involve three distinguished term-variables, usually denoted $S$, $P$ and $M$, in such a
way that $M$ occurs in both the diagrams in the
premises and does not in the conclusion, whereas $S$ and $P$
occur in the conclusion as well as in the premises. Verifying the validity of a syllogism consists in calculating the composite diagram of the concatenation of its premises, if these form a composable pair, and compare it with the diagram for the conclusion.\\
For example, the verification of the validity of the syllogism. $\aa{P}{M},\ee{S}{M}\vdash\ee{S}{P}$ is suggestively represented by the drawing

\begin{eqnarray}\label{sylloline}
\AxiomC{\xymatrix@R=1.5ex{S\ar[r]&\b&M\ar[l]&P\ar[l]}}
\UnaryInfC{\EE{S}{P}}
\DisplayProof
\end{eqnarray}
whereas the invalidity of the syllogism $\oo{P}{M},\ee{M}{S}\vdash\ii{S}{P}$ is confirmed by the fact that the pair $(P\leftarrow\b\arr\b\leftarrow M,M\arr\b\leftarrow S)$
although concatenable is not composable.

\begin{remark}\label{doppo}
Anticipating~\ref{rules} and~\ref{propopropo}, we haste to remark that in calculating the composite of a composable pair of diagrams no bullet symbol is deleted, so that the composite contains as many bullets as in the concatenation of the diagrams in the given pair. It
is useful, when one also takes into account the orientation of the involved arrow symbols, for rejecting an invalid form of syllogism, which can be 
rejected with the linear diagrams in~\cite{MR1149957} as well, but one has to go through all the $232$ invalid moods, as explained there. For instance, the syllogism $\oo{P}{M},\ee{M}{S}\vdash\ii{S}{P}$ is invalid since a single bullet symbol
occurs in the conclusion, whereas three of them occur in the premises. The syllogism $\aa{P}{M},\ii{S}{M}\vdash \ee{S}{P}$ is invalid since the syllogistic diagram for the conclusion contains a single bullet and a pair of arrows converging to it, whereas a single bullet and a pair of arrows diverging from it are contained in the syllogistic diagram for the second premise.
\end{remark}

For every term-variable $A$, particularly interesting
instances of syllogistic diagrams are the following:
\[\begin{tabular}{llllll}
\aa{A}{A}: & \AA{A}{A} && \ee{A}{A}: & \EE{A}{A}\\
\\
\ii{A}{A}: & \II{A}{A} && \oo{A}{A}: & \OO{A}{A}
\end{tabular}\]
where $\catebf{A}_{AA}$ and $\catebf{I}_{AA}$ are
referred to as \emph{laws of identity}, used by Aristotle without any explicit mention,
see~\cite{Lukasiewicz}. As will soonely be clear, the diagram for 
$\catebf{I}_{AA}$ represents existential import. 
The diagram for $\catebf{O}_{AA}$ is an expression of
the \emph{principle of contradiction}, which fact will be more
clearly illustrated in section~\ref{sec5}. The diagram for $\ee{A}{A}$ represents the emptyness of $A$.

\begin{definition}\label{rules}
The rules of inference of SYLL are the following:
\[\begin{tabular}{ccccc}
\AxiomC{\AA{A}{B}}\doubleLine\UnaryInfC{\Ao{A}{B}}\DisplayProof
&&\AxiomC{\EE{A}{B}}\UnaryInfC{\Eo{A}{B}}\DisplayProof&
\end{tabular}\]
\[\begin{tabular}{ccccc}
\AxiomC{\II{A}{B}}\UnaryInfC{\Io{A}{B}}\DisplayProof
&&\AxiomC{\OO{A}{B}}\doubleLine\UnaryInfC{\Oo{A}{B}}\DisplayProof&
\end{tabular}\]

\[\begin{tabular}{ccccc}
\AxiomC{\cate{D}A}
\AxiomC{A\cate{E}}
\BinaryInfC{\cate{D}A\cate{E}}
\DisplayProof&&
\AxiomC{A\cate{E}}
\AxiomC{\cate{D}A}
\BinaryInfC{\cate{D}A\cate{E}}
\DisplayProof&
\end{tabular}\]

\[\begin{tabular}{ccccc}
\AxiomC{$\cate{D}\arr A\arr\cate{E}$}
\UnaryInfC{$\cate{D}\arr\cate{E}$}
\DisplayProof&&
\AxiomC{$\cate{D}\leftarrow A\leftarrow\cate{E}$}
\UnaryInfC{$\cate{D}\leftarrow\cate{E}$}
\DisplayProof&
\end{tabular}\]
where the double line means that the rule can be used top-down as well as bottom-up. 
A \emph{proof tree} of SYLL is a tree where each node is a diagram and each branching is an instance of a rule of inference. 
A \emph{formal proof} of a syllogism is a proof tree
with its conclusion as the root and with each of its premises as leaves. A syllogism is \emph{provable} in SYLL if there is a formal proof for it.
\end{definition}

\begin{remark}
The last four rules in the previous definition can be equivalently substituted by the following:
\[\begin{tabular}{ccccc}
\AxiomC{\cate{D}\arr A}
\AxiomC{A\arr\cate{E}}
\BinaryInfC{\cate{D}\arr\cate{E}}
\DisplayProof&&
\AxiomC{A\arr\cate{E}}
\AxiomC{\cate{D}\arr A}
\BinaryInfC{\cate{D}\arr\cate{E}}
\DisplayProof&
\end{tabular}\]
\[\begin{tabular}{ccccc}
\AxiomC{$\cate{D}\leftarrow A$}
\AxiomC{$A\leftarrow\cate{E}$}
\BinaryInfC{$\cate{D}\leftarrow\cate{E}$}
\DisplayProof&&
\AxiomC{$A\leftarrow\cate{E}$}
\AxiomC{$\cate{D}\leftarrow A$}
\BinaryInfC{$\cate{D}\leftarrow\cate{E}$}
\DisplayProof&
\end{tabular}\]
\end{remark}

\begin{example}
The syllogism $\aa{P}{M},\ee{S}{M}\vdash\ee{S}{P}$ is provable, since 
a formal proof of it is
\[\AxiomC{\EE{S}{M}}
\AxiomC{\AA{P}{M}}
\UnaryInfC{\Ao{P}{M}}
\BinaryInfC{$\xymatrix{S\ar[r]&\b&M\ar[l]&P\ar[l]}$}
\UnaryInfC{\EE{S}{P}}
\DisplayProof\]
A different proof of the same syllogism is
\[\AxiomC{\AA{P}{M}}
\UnaryInfC{\Ao{P}{M}}
\AxiomC{\EE{S}{M}}
\BinaryInfC{$\xymatrix{S\ar[r]&\b&M\ar[l]&P\ar[l]}$}
\UnaryInfC{\EE{S}{P}}
\DisplayProof\]
\end{example}
\begin{notation}
Proof trees will be also written in line 
by forgetting some inessential pieces of information.
The proof tree of a syllogism $P_1,P_2\vdash C$, will be written as $(P_1)\sharp(P_2)\vdash(C)$. 
Drawings like~(\ref{sylloline}) will be formally considered as abbreviations of proof trees that we will henceforth freely use without any further comment. 
\end{notation}

\begin{lemma}\label{useful}
The composite of a composable pair whose components are syllogistic diagrams or reversals of them, is a syllogistic diagram in exactly the following cases:
\begin{itemize}
\item[(i)] $(\AA{S}{M},\AA{M}{P})$
\item[(ii)] $(\EE{S}{M},\Ao{P}{M})$
\item[(iii)] $(\AA{S}{M},\EE{M}{P})$
\item[(iv)] $(\Ao{M}{S},\II{M}{P})$
\item[(v)] $(\II{S}{M},\AA{M}{P})$
\item[(vi)] $(\II{S}{M},\EE{M}{P})$
\item[(vii)] $(\Ao{M}{S},\OO{M}{P})$
\item[(viii)] $(\OO{S}{M},\Ao{P}{M})$
\end{itemize}
\end{lemma}
\begin{proof}
Clearly, each of the listed composable pairs yields a composite syllogistic diagram involving only $S$ and $P$. Conversely, by also keeping in mind remark~\ref{doppo}, we proceed by cases:
\begin{itemize}
\item[(a)] the only way to obtain $S\arr P$ as a composite 
is by (i), since no bullet symbol is allowed to occur.
\item[(b)] the only way to obtain $S\arr\b\leftarrow P$ as a composite 
  is by either (ii) or (iii), since exactly one bullet symbol must occur with two
  arrow symbols converging to it.
\item[(c)] the only way to obtain $S\leftarrow\b\arr P$ as a
  composite is by either
  (iv) or (v), since exactly one bullet symbol must occur with
  two arrow symbols diverging from it.
\item[(d)] the only way to obtain $S\leftarrow\b\arr\b\leftarrow P$ as
  a composite is by either
  (vi), (vii) or (viii), since exactly two bullet symbols must occur 
  together with three alternating arrow symbols.
\end{itemize}
\end{proof}

\begin{remark}\label{rulesyl}
It is an easy exercise to read off the well-known \emph{rules of the syllogism}
from the list in lemma~\ref{useful}, also by taking into account remark~\ref{doppo}. 
\begin{itemize}
\item[(1)] From two negative premises nothing can be inferred.
\item[(2)] From two particular premises nothing can be inferred.
\item[(3)] If the first premise of a syllogism is particular, whereas its second premise is negative, then nothing can be inferred.
\item[(4)] If one premise is particular, then the conclusion is particular.
\item[(5)] The conclusion of a syllogism is negative if and only if so is one of its premises.
\end{itemize}
\end{remark}

The next theorem shows that the syllogisms in table~(\ref{questa}) are exactly those that are provable. The proof is purely syntactical and based on lemma~\ref{useful}. On one hand we proceed top-down constructing
a scheme of formal proof for any syllogism, from the syllogistic diagrams for its premises.
On the other hand we proceed bottom-up by cases, showing that the provable syllogisms leading to a possible syllogistic conclusion are among those of table~(\ref{questa}).
\begin{theorem}\label{propopropo}
A syllogism is valid if and only if it is provable in SYLL.
\end{theorem}
\begin{proof}
The syllogistic diagrams for the premises of a syllogism in table~(\ref{questa}), or their reversals, form composable pairs $(S\cate{A}\arr M, M\arr\cate{B}P)$ or $(S\cate{A}\leftarrow M,M\leftarrow\cate{B}P)$ that occurr among the ones listed in lemma~\ref{useful} and viceversa. Lemma~\ref{useful} ensures that the roots of the formal proofs
\[\AxiomC{$S\cate{A}\arr M$}
\AxiomC{$M\arr\cate{B}P$}
\BinaryInfC{$S\cate{A}\arr M\arr\cate{B}P$}
\UnaryInfC{$S\cate{A}\arr\cate{B}P$}
\DisplayProof\qquad
\AxiomC{$S\cate{A}\leftarrow M$}
\AxiomC{$M\leftarrow\cate{B}P$}
\BinaryInfC{$S\cate{A}\leftarrow M\leftarrow\cate{B}P$}
\UnaryInfC{$S\cate{A}\leftarrow\cate{B}P$}
\DisplayProof\]
are the syllogistic diagrams for the conclusion of any syllogism in table~(\ref{questa}).\\

By lemma~\ref{useful} (i), the only way to obtain $\catebf{A}_{SP}$ as a conclusion of a formal proof
is abbreviated as 
\[\AxiomC{\xymatrix{S\ar[r]&M\ar[r]& P}}
\UnaryInfC{\xymatrix{S\ar[r]&P}}
\DisplayProof\]
which amounts to the proof-tree
$(\catebf{A}_{MP})\sharp(\catebf{A}_{SM})\vdash(\catebf{A}_{SP})$
validating the mood $\catebf{AAA}$ in the first figure.\\
By lemma~\ref{useful} (ii) and (iii),
the only ways to obtain $\catebf{E}_{SP}$ as a conclusion of a formal proof, are
abbreviated as 
\[\AxiomC{\xymatrix{S\ar[r]&\b&M\ar[l]&P\ar[l]}}
\UnaryInfC{\xymatrix{S\ar[r]&\b&&P\ar[ll]}}
\DisplayProof
\qquad
\AxiomC{\xymatrix{S\ar[r]&M\ar[r]&\b&P\ar[l]}}
\UnaryInfC{\xymatrix{S\ar[rr]&&\b&P\ar[l]}}
\DisplayProof\]
The leftmost can be read as either the proof tree
$(\catebf{A}_{PM})^{ }\sharp(\catebf{E}_{SM})\vdash(\catebf{E}_{SP})$
or the proof tree
$(\catebf{A}_{PM})^{ }\sharp(\catebf{E}_{MS})^{ }\vdash(\catebf{E}_{SP})$
which validate the mood $\catebf{AEE}$ in the
second and fourth figures, respectively.
The rightmost can be read as either the proof tree
$(\catebf{E}_{MP})\sharp(\catebf{A}_{SM})\vdash(\catebf{E}_{SP})$
or the proof tree
$(\catebf{E}_{PM})^{ }\sharp(\catebf{A}_{SM})\vdash(\catebf{E}_{SP})$
which validate the mood $\catebf{EAE}$ in the first and second
figures, respectively.\\
By lemma~\ref{useful} (iv) and (v),
the only ways to obtain $\catebf{I}_{SP}$ as a conclusion of a formal proof
are abbreviated as 
\[\AxiomC{\xymatrix{S&M\ar[l]&\b\ar[l]\ar[r]&P}}
\UnaryInfC{\xymatrix{S&&\b\ar[ll]\ar[r]&P}}
\DisplayProof
\qquad
\AxiomC{\xymatrix{S&\b\ar[l]\ar[r]&M\ar[r]&P}}
\UnaryInfC{\xymatrix{S&\b\ar[l]\ar[rr]&&P}}
\DisplayProof\]
The leftmost can be read as either the proof tree
$(\catebf{I}_{MP})\sharp(\catebf{A}_{MS})^{ }\vdash(\catebf{I}_{SP})$
or the proof tree $(\catebf{I}_{PM})^{ }\sharp(\catebf{A}_{MS})^{ }\vdash(\catebf{I}_{SP})$
which validate the mood \catebf{IAI} in the third and fourth figures, respectively.
The rightmost can be read as either the proof tree $(\catebf{A}_{MP})\sharp(\catebf{I}_{SM})\vdash(\catebf{I}_{SP})$
or the proof tree $(\catebf{A}_{MP})\sharp(\catebf{I}_{MS})^{ }\vdash(\catebf{I}_{SP})$
that validate the mood \catebf{AII} in the first and third figures, respectively.\\
By lemma~\ref{useful} (vi), (vii) and (viii),
the only ways to obtain $\catebf{O}_{SP}$ as a conclusion of a formal proof are 
abbreviated as 
\[\AxiomC{\xymatrix{S&\b\ar[l]\ar[r]&M\ar[r]&\b&P\ar[l]}}
\UnaryInfC{\xymatrix{S&\b\ar[l]\ar[rr]&&\b&P\ar[l]}}
\DisplayProof\]
\[\AxiomC{\xymatrix{S&M\ar[l]&\b\ar[l]\ar[r]&\b&P\ar[l]}}
\UnaryInfC{\xymatrix{S&&\b\ar[ll]\ar[r]&\b&P\ar[l]}}
\DisplayProof\]
\[\AxiomC{\xymatrix{S&\b\ar[l]\ar[r]&\b&M\ar[l]&P\ar[l]}}
\UnaryInfC{\xymatrix{S&\b\ar[l]\ar[r]&\b&&P\ar[ll]}}
\DisplayProof\]
The first can be read as any of the proof trees $(\catebf{E}_{MP})\sharp(\catebf{I}_{SM})\vdash(\catebf{O}_{SP})$,
$(\catebf{E}_{PM})^{ }\sharp(\catebf{I}_{SM})\vdash(\catebf{O}_{SP})$,
$(\catebf{E}_{MP})\sharp(\catebf{I}_{MS})^{ }\vdash(\catebf{O}_{SP})$,
$(\catebf{E}_{PM})^{ }\sharp(\catebf{I}_{MS})^{ }\vdash(\catebf{O}_{SP})$,
that validate the mood \catebf{EIO} in all the
figures. The second can be read as the proof tree $(\catebf{O}_{MP})\sharp(\catebf{A}_{MS})^{}\vdash(\catebf{O}_{SP})$ that validates the mood \catebf{OAO} in
the third figure. The third can be read as the proof tree
$(\catebf{A}_{PM})^{ }\sharp(\catebf{O}_{SM})\vdash(\catebf{O}_{SP})$
validating the mood \catebf{AOO} in the second figure.
\end{proof}

Next is the extension of theorem~\ref{propopropo} to the strengthened syllogisms.
\begin{definition}
Let $\textup{SYLL}^+$ denote the formal system which is obtained from SYLL by the addition of the rule
\[\AxiomC{}
\UnaryInfC{\II{A}{A}}
\DisplayProof\]
to the rules in definition~\ref{rules} and with suitably extended notions of proof tree, formal proof and provability.
\end{definition}

\begin{lemma}\label{useful2}
The composite of a composable triple whose components are syllogistic
diagrams, reversals of them or existential imports, is a syllogistic diagram in exactly the following cases:
\begin{itemize}
\item[(i)] $(\II{S}{S},\AA{S}{M},\AA{M}{P})$
\item[(ii)] $(\Ao{M}{S},\II{M}{M},\AA{M}{P})$
\item[(iii)] $(\Ao{M}{S},\Ao{P}{M},\II{P}{P})$
\item[(iv)] $(\II{S}{S},\AA{S}{M},\EE{M}{P})$
\item[(v)] $(\II{S}{S},\EE{S}{M},\Ao{P}{M})$
\item[(vi)] $(\Ao{M}{S},\II{M}{M},\EE{M}{P})$
\end{itemize}
\end{lemma}
\begin{proof}
On one hand, it is clear that each of the listed composable triples
yields a syllogistic diagram as a composite. On the other hand,
by also keeping in mind remark~\ref{doppo}, we proceed by cases:
\begin{itemize}
\item[(a)] there is no way to obtain $S\arr P$ as the
  composite of a composable triple as in the statement, because of the
  occurrence of one indelible bullet symbol in any existential import
  for $S$, $M$ or $P$.
\item[(b)] there is no way to obtain $S\arr\b\leftarrow P$ as the
  composite of composable triple as in the statement, because of the presence of
  one indelible bullet symbol in any existential import for $S$, $M$ or
  $P$ together with two arrow symbols diverging from it.
\item[(c)] the only ways to obtain $S\leftarrow\b\arr P$ as the
  composite of composable triple as in the statement, under an
  existential import for $S$, $M$ or $P$, is by either (i), (ii) or (iii), since
exactly one bullet symbol must occur in the composite together with two
morphisms diverging from it.
\item[(d)] there is no way to obtain $S\leftarrow\b\arr\b\leftarrow P$
  as the composite of a composable triple as in the statement under an
  existential import for $P$, since such a composite would be of the form $S\cate{D}\leftarrow\b\arr P$ which by no means can be $S\leftarrow\b\arr\b\leftarrow P$. The only ways to obtain $S\leftarrow\b\arr\b\leftarrow P$ as a composite, 
under an existential import for $S$ or $M$ is by either (iv), (v) or (vi), since
exactly two bullet symbols must occur in the composite, together with three
  alternating morphisms.
\end{itemize}
\end{proof}
We end this section with the theorem which is the extension of 
theorem~\ref{propopropo} to the strengthened syllogisms. Its proof is completely analogous to the previous and is left to the reader, who is invited to carry it out on the base of lemma~\ref{useful2}.
\begin{theorem}\label{propos}
A strengthened syllogism is valid if and only if it is provable in $\textup{SYLL}^+$.
\end{theorem}

\section{On the Square of Opposition}\label{sec5}
We here want to point out the existing connections
between the so far described calculus and the laws of
the \emph{square of opposition}
\[\xymatrix{\aa{X}{Y}\ar@{--}[dddrrr]|{\textrm{contradiction}}
\ar@{-->}[ddd]_{\textrm{subalternation}}\ar@{--}[rrr]^{\textrm{contrariety}}&&&
\ee{X}{Y}\ar@{-->}[ddd]^{\textrm{subalternation}}\\
\\
\\
\ii{X}{Y}\ar@{--}[uuurrr]|{\hole}
\ar@{--}[rrr]_{\textrm{subcontrariety}}&&&\oo{X}{Y}}\]
in which 
\begin{itemize}
\item[-] \aa{X}{Y} and \oo{X}{Y}, as well as \ee{X}{Y} and
\ii{X}{Y}, are \emph{contradictory} because they negate each other and
in turn cannot hold together.
\item[-] under existential import, \aa{X}{Y} and \ii{X}{Y} as well as \ee{X}{Y} and \oo{X}{Y}, are \emph{subaltern} because \ii{X}{Y} is provable from
\aa{X}{Y} and \oo{X}{Y} is provable from \ee{X}{Y}, but not the converse,
  in both cases.
\item[-] under existential import \aa{X}{Y} and \ee{X}{Y} are \emph{contraries} because
the negation of each of them is provable from the other, but not the converse.
\item[-] under existential import \ii{X}{Y} and \oo{X}{Y} are \emph{subcontraries} because
each of them is provable from the negation of the other, but not the converse.
\end{itemize}
The laws of contradiction are the logical consequences
$\aa{X}{Y},\oo{X}{Y}\vdash\oo{X}{X}$, $\ee{X}{Y},\ii{X}{Y}\vdash\oo{X}{X}$
and this is the reason why we look at \oo{X}{X} as expressing contradiction in SYLL.
The remaining laws are condensed into the logical consequences
$\aa{X}{Y},\ii{X}{X}\vdash\ii{X}{Y}$, $\ee{X}{Y},\ii{X}{X}\vdash\oo{X}{Y}$
since they immediately provide the laws of subalternation. They express the laws of contrariety because \ii{X}{Y} is the negation of \ee{X}{Y} and \oo{X}{Y} is the negation of \aa{X}{Y}. They express the laws of subcontrariety since \aa{X}{Y} is the negation of \oo{X}{Y} and \ee{X}{Y} is the negation of \ii{X}{Y}.
\begin{proposition}\label{prop38}
The laws of the square of opposition are provable in $\textup{SYLL}^+$.
\end{proposition}
\begin{proof}
The laws of contradiction are provable by the proof trees
$(\aa{X}{Y})^{ }\sharp(\oo{X}{Y})\vdash(\oo{X}{X})$,
$(\ee{X}{Y})^{ }\sharp(\ii{X}{Y})\vdash(\oo{X}{X})$.
The remaining laws correspond to the proof trees
$(\aa{X}{Y})\sharp(\ii{X}{X})\vdash(\ii{X}{Y})$,
$(\ee{X}{Y})\sharp(\ii{X}{X})\vdash(\oo{X}{Y})$.
Both the proofs cannot be reversed
since one bullet symbol occurs in \ii{X}{Y} and no bullet symbols
occur in \aa{X}{Y}, two bullet symbols occur in \oo{X}{Y} and one bullet 
symbol occurs in \ee{X}{Y}.
\end{proof}

\section{Further discussion}\label{disc}
In this section we discuss informally 
the idea behind the syllogistic diagrams,
specifically the meaning of the arrow and bullet symbols, and the possibility of extending the calculus to $n$-term syllogisms, $n\geq 1$.\\

Concerning the first topic, we
do not have a complete answer and think that the subject deserves an 
investigation. Anyway, we can say that the syllogistic diagrams firstly
came out of an attempt to represent diagrammatically the degrees of
agreement and disagreement among terms of the natural language, as
conveyed by the categorical propositions. The syntactic primitives
that constitute them were employed since the beginning and, in
our opinion, one of their features is that they support an
intensional interpretation of terms, namely as concepts, rather than an
extensional one, namely as classes of individuals, which fact could be appreciable, see~\cite{MR2674955} and~\cite{MR1672695} for example. On the other hand, the controversy on the Aristotelian theory of syllogism being extensional or intensional could be considered as a futile one, see~\cite{Lukasiewicz}. In recent times we went aware of De Morgan's paper~\cite{ADeMorgan}, in which the so called ``spicular notation'' for the syllogistics was introduced, but also see~\cite{MR2464674}. De Morgan's aim was, among others, to extend Aristotle's syllogistics to complemented terms. The complement of a term $A$ is the term that means non-$A$, which De Morgan was used to denote with the corresponding lower case letter $a$. The syntactic primitives of De Morgan's system
are the symbols $\texttt{)}$, $\texttt{(}$, $\textcdot$, 
together with countably many term-variables $A,a,B,b,C,c\ldots$. De Morgan lets a term-variable be enclosed by a parethesis, as in $A\texttt{)}$ or $\texttt{(}A$, to express universal
quantification, that is ``all $A$s'', whereas he lets a term-variable
be excluded by a parethesis, as in $\texttt{)}A$ or $A\texttt{(}$,
to mean particular quantification, namely ``some $A$s''. In modern jargon, a
term-variable is said to be distributed in the first case and undistributed
in the second. Furthermore, he lets an even number of dots, or none at
all, between parentheses, express affirmation or agreement of terms, whereas he lets an odd number of dots express negation or disagreement of terms.
The following are the fundamental categorical 
propositions how they appear in the spicular notation:
\[\begin{tabular}{lllll}
\aa{A}{B}: & $A\spica B$ && \ee{A}{B}: & $A\spice B$\\
\\
\ii{A}{B}: & $A\spici B$ && \oo{A}{B}: & $A\spico B$
\end{tabular}\]
which accordingly should now be read

\[\begin{tabular}{l}
\aa{A}{B}: All $A$s are some $B$s\\
\\
\ee{A}{B}: All $A$s are not all $B$s\\
\\
\ii{A}{B}: Some $A$s are some $B$s\\
\\
\oo{A}{B}: Some $A$s are not all $B$s
\end{tabular}\]

We don't want to go now into a detailed comparison between our system SYLL and De Morgan's, but rather to point out that a way to give meaning to the syntactic primitives of SYLL could be based on the observation that the possibility of making a distinction between a term being distributed or not, as well as between affirmative and negative modes of predication, is supported by our diagrammatic formalism too, together with the possibility of handling complements of terms. Indeed, in our formalism a term-variable $A$ should be considered as distributed if fitting in a part such as $A\arr$ or $\leftarrow A$, whereas it should be considered as undistributed if fitting in a part such as $A\leftarrow $ or $\arr
A$. A term-variable $A$ should be considered as occurring in negated form if 
fitting in a part such as  $A\arr\b$ or $\b\leftarrow A$, both of which may be abbreviated as $a$. Thus we observe in passing that in our opinion the giving of an explicit enconding of negated terms through the syntactic primitives of SYLL is one of its remarkable features. 
On the other hand, a term-variable $A$ is in positive form if fitting in a part such as $\b\arr A$ or $A\leftarrow\b$. For example, in the syllogistic diagram for \aa{A}{B} the term-variable $A$ is distributed whereas the term-variable $B$ is not, or in that for \ee{A}{B} both the term-variables are distributed and negated, and similarly for the remaining syllogistic diagrams. Moreover, the rereading of the syllogistic diagrams under this perspective is in line with what happens for the linear diagrams in~\cite{MR1149957} for what concerns the obversion of the categorical propositions. By obversion ``No $A$ is $B$'' is equivalent to ``Each $A$ is non-$B$'' whereas ``Some $A$ is not $B$'' is equivalent to ``Some $A$ is non-$B$'', which fact is clearly expressed by the appearance of the syllogistic diagrams for \ee{A}{B} and \oo{A}{B}. 
By the introduction of complemented terms in syllogistics, De Morgan was able to introduce four more categorical propositions and also to let the
particular and universal affirmative modes of predication be the
fundamental ones. We conjecture that SYLL supports such an extension too, through the introduction of four further corresponding syllogistic diagrams. Finally, we end this digression by mentioning that each syllogistic diagram can be also more naively conceived as an ``abstract copula'':
\begin{quote} 
[\ldots] a formal mode of joining two terms which carries no meaning, and obeys no law except such as is barely necessary to make the forms of inference follow. See~\cite{ADeMorgan}.
\end{quote}

We end by briefly discussing the possibility of extending the calculus to
$n$-term syllogisms. This seems to be a peculiarity of syllogistic reasoning with linear diagrams, as observed in~\cite{MR1149957}.
For every natural number $n$, $n\geq 1$, an $n$-term syllogism is a 
logical consequence $P_1,\ldots,P_{n-1}\vdash P_n$ in which all the $P_i$'s are categorical propositions such that for every $1\leq i\lt n-1$ the categorical propositions $P_i$ and $P_{i+1}$ have exactly one term-variable in common. Thus, an $n$-term syllogism involves exactly $n$ term-variables $A_1,\ldots,A_n$, with $A_1$ in $P_{n-1}$ and $A_n$ in $P_1$, which are the subject and predicate of $P_n$, respectively. The total number of valid $n$-term syllogisms is
$3n^2-n$, see~\cite{Meredith}, where such a formula was obtained 
by rejecting the invalid moods on the base of the traditional rules of the
syllogism. The same formula has been reobtained by direct calculation in~\cite{Smyth}
and~\cite{MR0107598}. We conjecture that our system allows 
the retrieving of this result and moreover the extension of
theorems~\ref{propopropo} and~\ref{propos} to the case of $n$-term
syllogisms, but leave the investigation of these topics to a subsequent paper.
For the time being, the description of the valid $n$-term syllogisms 
for $n=1,2$ follows. For $n=1$ there is exactly one figure, that is
$A_1A_1$ and only two valid moods for it, that is $\catebf{A}$ and $\catebf{I}$ so that, as observed in~\cite{Lukasiewicz} and~\cite{Meredith}, the only valid $1$-term syllogisms are
$\vdash\catebf{A}_{A_1A_1}$ and $\vdash\catebf{I}_{A_1A_1}$, that is the laws of
identity we hinted at in section~\ref{sec4}. For $n=2$ there are two figures, as shown in the table
\[\begin{tabular}{|l|c|c|}
\hline
& fig. 1 & fig. 2\\
\hline
premise & $A_1A_2$ & $A_1A_2$\\
\hline
conclusion & $A_1A_2$ & $A_2A_1$\\
\hline
\end{tabular}\]
and ten valid $2$-term syllogisms, six in the first figure and four
in the second, as follows:
\begin{description}
\item[figure 1:] $\catebf{A}_{A_1A_2}\vdash\catebf{A}_{A_1A_2}$,
  $\catebf{E}_{A_1A_2}\vdash\catebf{E}_{A_1A_2}$,
  $\catebf{I}_{A_1A_2}\vdash\catebf{I}_{A_1A_2}$,
  $\catebf{O}_{A_1A_2}\vdash\catebf{O}_{A_1A_2}$, 
and the \emph{laws of subalternation}
  $\catebf{A}_{A_1A_2}\vdash\catebf{I}_{A_1A_2}$, $\catebf{E}_{A_1A_2}\vdash\catebf{O}_{A_1A_2}$ which both hold under existential import on $A_1$.
\item[figure 2:] $\catebf{E}_{A_1A_2}\vdash\catebf{E}_{A_2A_1}$,
  $\catebf{I}_{A_1A_2}\vdash\catebf{I}_{A_2A_1}$ which are the
  \emph{laws of simple conversion}, and
  $\catebf{A}_{A_2A_1}\vdash\catebf{I}_{A_1A_2}$,
  $\catebf{E}_{A_2A_1}\vdash\catebf{O}_{A_1A_2}$
  which are the \emph{laws of conversion per accidens}, that hold under existential import on $A_2$ and $A_1$, respectively.
\end{description}
In order to retrieve the law of identity $\vdash\aa{A_1}{A_1}$ the rule 
\[\AxiomC{}
\UnaryInfC{\AA{A_1}{A_1}}
\DisplayProof\]
has to be added in definition~\ref{rules}. For $n=2$, the laws of subalternation have been already proved in the previous section and, excluding the laws of conversion per accidens, the remaining syllogisms are immediate. We prove the laws of conversion per accidens:
\begin{itemize}
\item[-] $\catebf{A}_{A_2A_1}\vdash\catebf{I}_{A_1A_2}$
\[\AxiomC{}
\UnaryInfC{\II{A_2}{A_2}}
\AxiomC{\AA{A_2}{A_1}}
\BinaryInfC{$\xymatrix{A_2&\b\ar[r]\ar[l]&A_2\ar[r]&A_1}$}
\UnaryInfC{\II{A_2}{A_1}}
\UnaryInfC{\II{A_1}{A_2}}
\DisplayProof\]
\item[-] $\catebf{E}_{A_2A_1}\vdash\catebf{O}_{A_1A_2}$
\[\AxiomC{\EE{A_2}{A_1}}
\AxiomC{}
\UnaryInfC{\II{A_1}{A_1}}
\BinaryInfC{$\xymatrix{A_2&\b&A_1\ar[l]&\b\ar[r]\ar[l]&A_1}$}
\UnaryInfC{\Oo{A_1}{A_2}}
\UnaryInfC{\OO{A_1}{A_2}}
\DisplayProof\]

\end{itemize}

\bibliographystyle{plain}
\bibliography{BiblioTeX}
\end{document}